\documentstyle[12pt]{article}
\input{amssym}

\newcommand{\p}{\rm \partial}
\begin{document}
\title{A symmetry classification  for a class of\\
$(2+1)$-nonlinear wave equation}
\author{M. Nadjafikhah\thanks{e-mail: m\_nadjafikhah@iust.ac.ir
Department of Mathematics, Iran University of Science and
Technology, Narmak, Tehran, IRAN.} \and R. Bakhshandeh
Chamazkoti\thanks{e-mail: r\_bakhshandeh@iust.ac.ir}\and A.
Mahdipour-Shirayeh\thanks{e-mail: mahdipour@iust.ac.ir}}
\date{}
\maketitle
\renewcommand{\sectionmark}[1]{}
\begin{abstract}
In this paper, a symmetry classification of a $(2+1)$-nonlinear
wave equation $u_{tt}-f(u)(u_{xx}+u_{yy})=0$ where $f(u)$ is a
smooth function on $u$, using Lie group method, is given. The
basic infinitesimal method for calculating symmetry groups is
presented, and used to determine the general symmetry group of
this $(2+1)$-nonlinear wave equation.
\end{abstract}
\section{Introduction}
It is well known that the symmetry group method plays an important
role in the analysis of differential equations.   The history of
group classification methods goes back to Sophus Lie. The first
paper on this subject is \cite{[1]}, where Lie proves that a
linear two-dimensional second-order PDE may admit at most a
three-parameter invariance group (apart from the trivial infinite
parameter symmetry group, which is due to linearity).  He computed
the maximal invariance group of the one-dimensional heat
conductivity equation and utilized this symmetry to construct its
explicit solutions. Saying it the modern way, he performed
symmetry reduction of the heat equation. Nowadays symmetry
reduction is one of the most powerful tools for solving nonlinear
partial differential equations (PDEs). Recently, there have been
several generalizations of the classical Lie group method for
symmetry reductions. Ovsiannikov \cite{[2]} developed the method
of partially invariant solutions. His approach is based on the
concept of an equivalence group, which is a Lie transformation
group acting in the extended space of independent variables,
functions and their derivatives, and preserving the class of
partial differential equations under study.
\newline $~~~~$ For many nonlinear
systems, there are only explicit exact solutions available. These
solutions play an important role in both mathematical analysis and
physical applications of the systems. There are a number of papers
to study $(1+1)$-nonlinear wave equations from the point of view
of Lie symmetries method. First, for solving some of the physical
problems, the quasi-linear hyperbolic equation with the form
\begin{eqnarray}
u_{tt}=[f(u)u_x]_x,\label{eq:1'}
\end{eqnarray}
in \cite{[3]} and later its the generalized cases
\begin{eqnarray}
u_{tt}=[f(x,u)u_x]_x,\hspace{1cm} u_{tt}=[f(u)u_x+g(x,u)]_x,
\end{eqnarray}
in \cite{[4]} and \cite{[5]}, respectively, are investigated. Also
the most important classes of the $(1+1)$-nonlinear wave equations
with the forms
\begin{eqnarray}
v_{tt}=f(x,v_x)v_{xx}+g(x,v_x),\hspace{1cm}
u_{tt}=f(x,u)u_{xx}+g(x,u),
\end{eqnarray}
can be found in two attempts \cite{[6]} and
 \cite{[7]} respectively. An alternative form of Eq. (\ref{eq:1'}) was
also investigated by Oron and Rosenau \cite{[8]} and Suhubi and
Bakkaloglu \cite{[9]}.  The equations
\begin{eqnarray}
u_{tt}=F(u)u_{xx},\hspace{1mm}u_{tt}+K(u)u_t=F(u)u_{xx},\hspace{1mm}u_{tt}+K(u)u_t=F(u)u_{xx}+H(u)u_x,
\end{eqnarray}
are  classified in \cite{[10],[11],[12]}, respectively. Lahno et
al. \cite{[13]} presented the most extensive list of symmetries of
the equations
\begin{eqnarray}
u_{tt}=u_{xx}+ F (t, x, u, u_x) ,
\end{eqnarray}
by using the infinitesimal Lie method, the technique of
equivalence transformations, and the theory of classification of
abstract low-dimensional Lie algebras. There are also some papers
\cite{[14],[15],[16]} devoted to the group classification of the
equations of the following form:
\begin{eqnarray}
u_{tt}=F(u_{xx}),\hspace{1cm} u_{tt} = F (u_x)u_{xx} +
H(u_x),\hspace{1cm} u_{tt} + u_{xx} = g(u, u_x),
\end{eqnarray}
Studies have also been made for $(2+1)$-nonlinear wave equation
with constant coefficients \cite{[17],[18],[19]}. In the special
case the $(2+1)$-dimensional nonlinear wave equation
\begin{equation}
u_{tt}=u^n(u_{xx}+u_{yy}),
\end{equation}
is investigated in \cite{[20]}. The goal of this paper is to
investigate the Lie symmetries for some class of $(2+1)$-nonlinear
wave equation
\begin{equation}
u_{tt}-f(u)(u_{xx}+u_{yy})=0,\label{eq:1}
\end{equation}
where $f(u)$ is a arbitrary smooth function of the variable $u$.
Clearly, in Eq. (\ref{eq:1}) case of $f_u=0$ namely
$f(u)=constant$ is not interest because this case reduces the wave
equation to a linear one. Similarly technics were applicable for
some classes of the nonlinear heat equations in \cite{[21],[22]}.
%
\section{ Symmetry Methods}
Let a partial differential equation contains  one dependent
variable and $p$ independent variables. The one-parameter Lie
group of transformations
\begin{eqnarray}
\overline{x}_i=x_i+\epsilon\xi_i(x,u)+O(\epsilon^2);\hspace{1cm}
\overline{u}=u+\epsilon\varphi(x,u)+O(\epsilon^2),\label{eq:3}
\end{eqnarray}
where $i=1,\ldots,p$, and
$\xi_i=\frac{{\p}\overline{x}_i}{{\p}\epsilon}\big|_{\epsilon=0}$,
acting on $(x,u)$-space has as its infinitesimal generator
\begin{eqnarray}
{\bf
v}=\xi_i\frac{\p}{{\p}x_i}+\varphi\frac{\p}{{\p}u},\hspace{1cm}i=1,\ldots,p~.\label{eq:4}
\end{eqnarray}
Therefore, the characteristic of the vector field ${\bf v}$ given
by (\ref{eq:4}) is the function
\begin{eqnarray}
Q(x,u^{(1)})=\varphi(x,u)-\sum_{i=1}^p\xi_i(x,u)\frac{{\p}u}{{\p}x_i}.\label{eq:5}
\end{eqnarray}
The symmetry generator associated with (\ref{eq:4}) given by
\begin{eqnarray}
{\bf v}=\xi\frac{\p}{{\p}x}+\eta\frac{\p}{{\p}y}+
\tau\frac{\p}{{\p}t}+\varphi\frac{\p}{{\p}u}~.\label{eq:6}
\end{eqnarray}
The second prolongation of ${\bf v}$ is the vector field
\begin{eqnarray}
{\bf v}^{(2)}={\bf
v}+\varphi^x\frac{\p}{{\p}u_x}+\varphi^y\frac{\p}{{\p}u_y}+\varphi^t\frac{\p}{{\p}u_t}+
\varphi^{xx}\frac{\p}{{\p}u_{xx}}+\varphi^{xy}\frac{\p}{{\p}u_{xy}}
+\varphi^{xt}\frac{\p}{{\p}u_{xt}}\\\nonumber
+\varphi^{yy}\frac{\p}{{\p}u_{yy}}
+\varphi^{yt}\frac{\p}{{\p}u_{yt}}+\varphi^{tt}\frac{\p}{{\p}u_{tt}}~.\label{eq:7}
\end{eqnarray}
with coefficient
\begin{eqnarray}
\varphi^{\iota}&=&D_{\iota}Q+\xi u_{x\iota}+\eta u_{y\iota}+\tau
u_{t\iota},\label{eq:8}\\
\varphi^{\iota\jmath}&=&D_{\iota}(D_{\jmath}Q) +\xi
u_{x\iota\jmath}+\eta u_{y\iota\jmath}+\tau
u_{t\iota\jmath},\label{eq:9}
\end{eqnarray}
where $Q=\varphi-\xi u_x-\eta u_y-\tau u_t$ is the characteristic
of the vector field ${\bf v}$ given by (\ref{eq:6}) and $D_i$
represents total derivative and subscripts of $u$ are derivative
with respect to the respective coordinates. $\iota$ and $\jmath$
in above could be $x,y$ or $t$ coordinates. By the theorem 6.5. in
\cite{[23]}, ${\bf v}^{(2)}[u_{tt}-f(u)(u_{xx}+u_{yy})]=0$
whenever
\begin{eqnarray}
u_{tt}-f(u)(u_{xx}+u_{yy})=0.\label{eq:10}
\end{eqnarray}
Since $${\bf
v}^{(2)}[u_{tt}-f(u)(u_{xx}+u_{yy})]=\varphi^{tt}-\varphi
f_u(u_{xx}+u_{yy}) -f(u)(\varphi^{xx}+\varphi^{yy}),$$ therefore
\begin{eqnarray}
\varphi^{tt}-\varphi f_u(u_{xx}+u_{yy})
-f(u)(\varphi^{xx}+\varphi^{yy})=0.\label{eq:11}
\end{eqnarray}
Using the formula (\ref{eq:9}) we obtain coefficient functions
$\varphi^{xx}, \varphi^{yy}, \varphi^{tt}$ as
\begin{eqnarray}
\varphi^{xx}&=&D_{x}^2Q +\xi u_{xx}+\eta u_{yxx}+\tau
u_{txx},\label{eq:12}\\
\varphi^{yy}&=&D_{y}^2Q +\xi u_{xyy}+\eta u_{yyy}+\tau
u_{tyy},\label{eq:9'}\\
\varphi^{tt}&=&D_{t}^2Q +\xi u_{xtt}+\eta u_{ytt}+\tau
u_{ttt},\label{eq:9''}
\end{eqnarray}
where the operators $D_x$, $D_y$ and $D_t$ denote the total
derivatives with respect to $x,y$ and $t$:
\begin{eqnarray}\nonumber
D_x&=&\frac{\p}{{\p}x}+u_x\frac{\p}{{\p}u}+u_{xx}\frac{\p}{{\p}u_x}+u_{xy}\frac{\p}{{\p}u_y}+
u_{xt}\frac{\p}{{\p}u_t}+\ldots\\
D_y&=&\frac{\p}{{\p}y}+u_y\frac{\p}{{\p}u}+u_{yy}\frac{\p}{{\p}u_y}+u_{yx}\frac{\p}{{\p}u_x}+
u_{yt}\frac{\p}{{\p}u_t}+\ldots\\
D_t&=&\frac{\p}{{\p}t}+u_t\frac{\p}{{\p}u}+u_{tt}\frac{\p}{{\p}u_t}+u_{tx}\frac{\p}{{\p}u_x}+
u_{ty}\frac{\p}{{\p}u_y}+\ldots\nonumber
\end{eqnarray}
and by substituting them into invariance condition (\ref{eq:11}),
we are left with a polynomial equation involving the various
derivatives of $u(x,y,t)$ whose coefficients are certain
derivatives of $\xi,\eta,\tau$ and $\varphi$. Since
$\xi,\eta,\tau,\varphi$ only depend on $x,y,t,u$ we can equate the
individual coefficients to zero, leading to the complete set of
determining equations:
\begin{eqnarray}
\xi&=&\xi(x,t)\label{eq:15}\\
\eta&=&\eta(y,t)\label{eq:16}\\
\tau&=&\tau(x,y,t)\label{eq:17}\\
\varphi&=&\alpha(x,y,t)u+\beta(x,y,t)\label{eq:18}\\
\tau_t&=&\varphi_u=\alpha(x,y,t),\label{eq:18'}\\
\xi_{tt}&=&f(u)(\xi_{xx}-2\varphi_{xu})\label{eq:19}\\
\eta_{tt}&=&f(u)(\eta_{yy}-2\varphi_{yu})\label{eq:20}\\
\tau_{tt}&=&f(u)(\tau_{xx}+\tau_{yy})+2\varphi_{tu}\label{eq:21}\\
f_u\varphi&=&2f(u)(\xi_x-\tau_t)\label{eq:22}\\
f_u\varphi&=&2f(u)(\eta_y-\tau_t)\label{eq:23}\\
f(u)\tau_x&=&\xi_t\label{eq:24}\\
f(u)\tau_y&=&\eta_t\label{eq:25}\\
\varphi_{tt}&=&f(u)(\varphi_{xx}+\varphi_{yy})\label{eq:26}
\end{eqnarray}
\section{classification of symmetries of the
model}
In this section we start to classify the symmetries of the
nonlinear wave equation (\ref{eq:1}). To fined a complete solution
of the above system we consider Eq. (\ref{eq:22}) and with
assumption $f_u\neq0$ we rewrite:
\begin{eqnarray}
\varphi=2\frac{f}{f_u}(\xi_x-\tau_t)\label{eq:27}
\end{eqnarray}
Note the case of $f(u)=constant$ explained in introduction. Two
general cases are possible:
\begin{eqnarray}
{\rm i})~~~~\frac{f}{f_u}&=&c,\label{eq:28}\\
{\rm ii})~~~~ \frac{f}{f_u}&=&g(u),\label{eq:29}
\end{eqnarray}
where $c$ is a constant.
\subsection{Case {\rm (i)}}
In this case with integrating from Eq. (\ref{eq:28}) with respect
to $u$ to obtain
\begin{eqnarray}
f(u)=Ke^{\frac{u}{c}},\label{eq:30}
\end{eqnarray}
where $K$ is an integration constant. Then the Eq. (\ref{eq:27})
reduce to
\begin{eqnarray}
\varphi=2c(\xi_x-\tau_t).\label{eq:31}
\end{eqnarray}
With substituting (\ref{eq:31}) into (\ref{eq:18'})-(\ref{eq:25})
we have
\begin{eqnarray}
&\xi(x)=c_1x+c_2; & \hspace{1cm}\eta(y)=c_1y+c_3;\\\nonumber
&\tau(t)=c_4t+c_5; &\hspace{1cm}\varphi=2c(c_1-c_4).\label{eq:32}
\end{eqnarray}
where $c_i$, $i=1,\ldots,5$, are arbitrary constants. The Lie
symmetry generator for Eq. (\ref{eq:1}) in this case (i) is
\begin{eqnarray}
{\bf v}=(c_1x+c_2)\frac{\p}{{\p}x}+(c_1y+c_3)\frac{\p}{{\p}y}+
(c_4t+c_5)\frac{\p}{{\p}t}+2c(c_1-c_4)\frac{\p}{{\p}u}~.\label{eq:33}
\end{eqnarray}
Therefore the symmetry algebra of the $(2+1)$-nonlinear wave
equation (\ref{eq:1}) is spanned by the vector fields
\begin{eqnarray}
{\bf v}_1=x{\p}_x+y{\p}_y+2c{\p}_u; & \hspace{1mm}{\bf
v}_2={\p}_x;\hspace{1mm}{\bf v}_3={\p}_y; \hspace{1mm}{\bf
v}_4=t{\p}_t-2c{\p}_u; \hspace{1mm}{\bf v}_5={\p}_t.\label{eq:34}
\end{eqnarray}
The commutation relations satisfied by generators (\ref{eq:34}) in
the case (i) are shown in table 1.
\begin{table}
\caption{Commutation relations satisfied by infinisimal generators
in Cases (i) and (ii)}\label{tab1}
\begin{tabular}{llllll} \hline
$[{\bf v}_i,{\bf v}_j]$ & ${\bf v}_1$ & ${\bf v}_2$& ${\bf v}_3$& ${\bf v}_4$& ${\bf v}_5$\\
\hline
${\bf v}_1$ & $0$          & $-{\bf v}_2$& $-{\bf v}_3$ & $0$         & $0$           \\
${\bf v}_2$ & ${\bf v}_2$  & $0$         & $0$          & $0$         & $0$           \\
${\bf v}_3$ & ${\bf v}_3$  & $0$         & $0$          & $0$         & $0$           \\
${\bf v}_4$ & $0$          & $0$         & $0$          & $0$         & $-{\bf v}_5$  \\
${\bf v}_5$ & $0$          & $0$         & $0$          & ${\bf v}_5$ & $0$            \\
\hline
\end{tabular}
\end{table}
The invariants associated with the infinitesimal generator ${\bf
v}_1$ are obtained by integrating the characteristic equation:
\begin{eqnarray}
\frac{dx}{x}=\frac{dy}{y}=\frac{dt}{0}=\frac{du}{2c}~.\label{eq:35}
\end{eqnarray}
and have the forms
\begin{eqnarray}
r=\frac{y}{x},~~ s=t,~~{\rm and}~~ \omega(r,s)=u(x,y,t)-2c\ln
x,\label{eq:36}
\end{eqnarray}
With substituting (\ref{eq:36}) into (\ref{eq:10}) to determine
the form of the function $\omega$ to obtain
\begin{eqnarray}
\omega_{ss}=Ke^{\frac{\omega}{c}}\Big((1+r^2)\omega_{rr}+
2r\omega_r-2c\Big),\label{eq:37}
\end{eqnarray}
By solving this partional differential equation we obtain  the
reduced equation
\begin{eqnarray}
\omega(r,s)=\zeta_1(r)+\zeta_2(s),\label{eq:38}
\end{eqnarray}
where $\zeta_1$ and $\zeta_2$ satisfy in following second-order
differential equations
\begin{eqnarray}
\ddot{\zeta_1}(r^2+1)+c_1e^{-{\zeta_1\over
c}}+2(r\dot{\zeta_1}-c)=0
;\hspace{1cm}\ddot{\zeta_2}+Kc_1e^{\zeta_2\over c}=0,\label{eq:39}
\end{eqnarray}
with $c_1, c, K,$ arbitrary constants. The characteristic equation
associated with ${\bf v}_4$ is
\begin{eqnarray}
\frac{dx}{0}=\frac{dy}{0}=\frac{dt}{t}=\frac{du}{-2c}~,\label{eq:40}
\end{eqnarray}
which generate the invariants $x$,~$y$,~$t^{-2c}e^{-u}$. Then the
similarity solution is chosen to have the form
\begin{eqnarray}
u(x,y,t)=2c\ln\frac{h(x,y)}{t}~.\label{eq:41}
\end{eqnarray}
By substituting (\ref{eq:41}) into (\ref{eq:10}) to determine the
form of the function $h$ to obtain
\begin{eqnarray}
{1\over K}-h(x,y)(h_{xx}+h_{yy})+h_x^2+h_y^2=0,\label{eq:42}
\end{eqnarray}
which has the solution
\begin{eqnarray}
h(x,y)=mx+py+q;~~~m^2+p^2=K^{-1},\label{eq:43}
\end{eqnarray}
where $m,~p,~q$ are arbitrary constants. For the remaining
infinitesimal generators ${\bf v}_2$, ${\bf v}_3$, ${\bf v}_5$,
the invariants associated are the arbitrary functions
$\lambda(y,t,u)$, $\mu(x,t,u)$, and $\nu(x,y,u)$ respectively.
\subsection{Case {\rm (ii)}}
In this case we classify solution of the wave equation
(\ref{eq:1}), with assumption $g_u\neq0$. With substituting
(\ref{eq:26}) into (\ref{eq:19})-(\ref{eq:20}), since $\xi$,
$\eta$ and $\tau$ are not dependent to $u$, therefore from
\begin{eqnarray}
u(x,y,t,u)=2g(u)(\xi_x-\tau_t),\label{eq:44}
\end{eqnarray}
and also from (\ref{eq:1}) and (\ref{eq:18}), we conclude
\begin{eqnarray}
g(u)=e_1u+e_2,\label{eq:45}
\end{eqnarray}
where $e_1\neq0$ and $e_2$ are arbitrary constants. Now we
substitute (\ref{eq:45}) into  (\ref{eq:29}) and rewrite
\begin{eqnarray}
\frac{f_u}{f}=\frac{1}{e_1u+e_2}.\label{eq:46}
\end{eqnarray}
Therefore by integrating from (\ref{eq:46}) with respect to $u$ we
have
\begin{eqnarray}
f(u)=L(e_1u+e_2)^{\frac{1}{e_1}},\label{eq:47}
\end{eqnarray}
where $L$ is an integration constant. Now by considering Eq.
(\ref{eq:15})-(\ref{eq:26}), it is not hard to find that the
components $\xi$, $\eta$, $\tau$ and $\varphi$ of infinitesimal
generators become
\begin{eqnarray}
&\xi(x)=c_1x+c_2; & \hspace{1cm}\eta(y)=c_1y+c_3;\\\nonumber
&\tau(t)=c_4t+c_5;
&\hspace{1cm}\varphi=2e_1(c_1-c_4)u+2e_2(c_1-c_4),\label{eq:48}
\end{eqnarray}
where $c_i$, $i=1,\ldots,5$,  are arbitrary constants. From above
the five infinitesimal generators can be constructed:
\begin{eqnarray}\nonumber
&{\bf v}_1=x{\p}_x+y{\p}_y+(2e_1u+2e_2){\p}_u; & \hspace{1cm}{\bf
v}_2={\p}_x;  \hspace{1cm}{\bf v}_3={\p}_y; \\
\hspace{1cm}&{\bf v}_4=t{\p}_t-2(e_1u+e_2){\p}_u;
&\hspace{1cm}{\bf v}_5={\p}_t.\label{eq:49}
\end{eqnarray}
It is easy to check that the infinitesimal generators
(\ref{eq:49}) from a closed Lie algebra whose it's corresponding
commutation relations are coincided with obtained results in table
1. For generator ${\bf v}_1$, the associated equations are
\begin{eqnarray}
\frac{dx}{x}=\frac{dy}{y}=\frac{dt}{0}=\frac{du}{2(e_1u+e_2)}~,\label{eq:50}
\end{eqnarray}
which generate the invariants $\displaystyle{p=\frac{y}{x}}$,
$q=t$, and $\displaystyle{\vartheta(p,q)=(u+{e_2\over
e_1})x^{-2e_1}}$. Consequently, the similarity solution is chosen
to have the form
\begin{eqnarray}
u(x,y,t)=\vartheta(t,{y\over
x})x^{2e_1}+\frac{e_2}{e_1}~.\label{eq:51}
\end{eqnarray}
We substitute (\ref{eq:51}) into (\ref{eq:10}) to obtain following
partional differential equation
\begin{eqnarray}
\vartheta_{pp}=L\vartheta^{1\over
e_1}[(q^2+1)\vartheta_{qq}+2q\vartheta_q(1-2e_1)+2e_1(2e_1-1)\vartheta],\label{eq:51}
\end{eqnarray}
as an example, for particular case $e_1 = 1$, the solution of
(\ref{eq:51}) is
\begin{eqnarray}
\vartheta(p,q)=\varsigma_1(p)\cdot\varsigma_2(q)~,\label{eq:52}
\end{eqnarray}
where $\varsigma_1(p)$ and $\varsigma_2(q)$ satisfy in second
order equations
\begin{eqnarray}
\ddot{\varsigma_1}-c\varsigma_1^2=0;~~~(q^2+1)\ddot{\varsigma_2}-
2q\dot{\varsigma_2}+2\varsigma_2-cL^{-1}=0,\label{eq:53}
\end{eqnarray}
where $c$ is a  arbitrary constant. Also characteristic equation
corresponding generator ${\bf v}_4$ is
\begin{eqnarray}
\frac{dx}{0}=\frac{dy}{0}=\frac{dt}{t}=\frac{du}{-2(e_1u+e_2)}~,\label{eq:54}
\end{eqnarray}
and so
\begin{eqnarray}
u(x,y,t)=l(x,y)t^{-2e_1}-\frac{e_2}{e_1}~.\label{eq:55}
\end{eqnarray}
Substitute (\ref{eq:55}) into (\ref{eq:10}), $l(x,y)$ satisfies in
the following equation:
\begin{eqnarray}
L(l_{xx}+l_{yy})l^{(e_1^{-1}-1)}-2e_1(2e_1+1)=0~.\label{eq:56}
\end{eqnarray}
For the remaining infinitesimal generators ${\bf v}_2$, ${\bf
v}_3$, ${\bf v}_5$, the invariants associated are the arbitrary
functions $r(y,t,u)$, $m(x,t,u)$, and $n(x,y,u)$ respectively.
\section{Conclusion and new ideas}
In this paper we have obtained some particular Lie point
symmetries group of the $(2+1)$-nonlinear wave equation
$u_{tt}-f(u)(u_{xx}+u_{yy})=0$ where $f(u)$ is a smooth function
on $u$, by using here the classical Lie symmetric method. In
section 2, the complete set of determining equations was obtained
by substituting the equations (\ref{eq:12}), (\ref{eq:9'}) and
(\ref{eq:9''}) in invariance condition (\ref{eq:11}) and then in
section 3, we classify the symmetries of this nonlinear wave
equation by assumption two cases in (\ref{eq:28}) and
(\ref{eq:29}) to consider $\displaystyle{\frac{f}{f_u}}$  is a
constant or is a smooth function with respect to $u$ and
$f_u\neq0$. The commutation relations satisfied by infinitesimal
generators in two cases are given in table 1, and their invariants
associated with the infinitesimal generators  are obtained. This
method is suitable for preliminary group classification of some
class of nonlinear wave equations \cite{[6],[7]}.
\newline $~~~~$There are some classes of $(2 +1)$-nonlinear wave equations that will be
investigated  by both classical or nonclassical  symmetries method
similarly whose we do for classical case. For examples
\begin{eqnarray}
u_{tt}-f(x,u)(u_{xx}+u_{yy})=0,\\\label{eq:56}
u_{tt}-f(x,u_x)(u_{xx}+u_{yy})=0,\label{eq:57}
\end{eqnarray}
or generalized case
\begin{eqnarray}
u_{tt}-f(x,y,u,u_x)(u_{xx}+u_{yy})=0,\label{eq:58}
\end{eqnarray}
are interested.
\section{Acknowledgements}
It is a pleasure to thank the anonymous referees for their
constructive suggestions and helpful comments which have
materially improved the presentation of the paper.
%

\end{document}